\documentclass[12pt]{amsart}

\usepackage{amsthm}

\usepackage{amssymb}
\usepackage{verbatim}
\usepackage[curve,matrix,arrow]{xy}
\usepackage{amsmath,amsfonts}
\usepackage{graphicx}
\usepackage{epsf}

\usepackage{amssymb, graphics}

\newcommand{\mathsym}[1]{{}}

\newcommand{\thmref}[1]{Theorem~\ref{#1}}

\newcommand{\eqnref}[1]{Equation~(\ref{#1})}
\newcommand{\remref}[1]{Remark~\ref{#1}}

\newcommand{\figref}[1]{Figure~\ref{#1}}
\newcommand{\tabref}[1]{Table~\ref{#1}}

  {\end{list}}

\def\ZZ{{\mathbb Z}}

\def\li{L_{i}}

\def\NN{{\mathbb N}}
\def\RR{{\mathbb R}}
\def\ZZ{{\mathbb Z}}

\newtheorem{theorem}{Theorem}[section]

\newtheorem{lemma}[theorem]{Lemma}

\theoremstyle{example}

\newtheorem{remark}[theorem]{Remark}

\theoremstyle{definition}

\theoremstyle{notation}

\newcommand{\dd}[1]{\delta_{#1}}

\newcommand{\nn}[1]{\mu_{#1}}

\newcommand{\ga}{\Gamma}

\newcommand{\tg}{\tau(\Gamma)}

\newcommand{\ee}[1]{E(#1)}
\newcommand{\vv}[1]{V(#1)}
\newcommand{\va}{\upsilon}

\newcommand{\gc}{g(C)}

\newcommand{\tc}{\theta}

\def\ZZ{{\mathbb Z}}

\def\<{\langle }
\def\>{\rangle }

\def\ed{\epsilon_{D}}

\def\elg{\ell (\ga)}

\def\gc{\bar{g}}
\def\ed{\epsilon(\ga)}

\def\bq{\textbf{q}}
\def\vg{\varphi (\ga)}
\def\elg{\ell (\ga)}
\newcommand{\tcg}{\theta (\ga)}
\def\lag{\lambda (\ga)}

\newcommand\T{\rule{0pt}{2.6ex}}
\newcommand\B{\rule[-1.2ex]{0pt}{0ex}}

\begin{document}

\title[Admissible Invariants of genus 3 Curves]
{Admissible Invariants of genus 3 Curves}

\author{Zubeyir Cinkir}
\address{Zubeyir Cinkir\\
Department of Mathematics\\
Zirve University\\
27260, Gaziantep, TURKEY\\}
\email{zubeyir.cinkir@zirve.edu.tr}



\begin{abstract}
Several invariants of polarized metrized graphs and their applications in Arithmetic Geometry are studied recently.
In this paper, we explicitly calculated these admissible invariants for all curves of genus $3$.
We find the sharp lower bound for the invariants $\varphi$, $\lambda$ and $\epsilon$ for all polarized metrized graphs of genus $3$.
This improves the lower bound given for Effective Bogomolov Conjecture for such curves.
\end{abstract}

\maketitle

\section{Introduction}\label{sec introduction}

Invariants of a polarized metrized graph $(\ga, \bq)$ are of interest for more than last twenty years because of their applications in arithmetic geometry and number theory. The invariants $\vg$, $\lag$ and $\ed$, see below for their definitions, are studied for their connection to the self intersection of admissible dualizing sheaf associated to a curve of genus at least $2$ over a global field. One can consult to S. Zhang's articles \cite{Zh1} and \cite{Zh2} for technical details. On the other hand, L. Szpiro showed in \cite{LS1} and \cite{LS2} that both Bogomolov and Effective Mordell conjectures can be stated as `suitable' upper and lower bounds to self intersection of  certain dualizing  sheaf.

Lower bounds for  $\vg$, $\lag$ are given in \cite[Theorems 2.11 and 2.13]{C5} for all curves of genus greater than $1$. However, obtaining sharp bounds and explicit computations for each possible case were done for only certain type of curves.
For curves of small genus, the invariants $\vg$, $\lag$ and $\ed$ are studied in \cite{AM2}, \cite{AM3}, \cite{KY1}, \cite{KY3}, \cite{J}, \cite{Fa}. It will be of reader's interest to see \cite{J2}, \cite{UK} and \cite{KY2}, too.

When $(\ga, \bq)$ has genus $2$, A. Moriwaki computed several invariants including $\ed$ for each case and obtained certain bound equivalent to the sharp lower bound $\vg \geq \frac{1}{27} \elg$, where $\elg$ is the total length of $\ga$. This bound is verified with different methods by X. Faber in \cite{Fa} and the author in \cite{C5}. For a such $\ga$, Moriwaki in \cite{AM3} and Jong in \cite{J} explicitly computed these invariants for each genus $2$ curves. Jong's work also extends to archimedean case, which is essential for number theoretic applications.

When $(\ga, \bq)$ has genus $3$, K. Yamaki in \cite{KY3}, Faber in \cite{Fa} and the author in \cite{C5} studied these invariants. Moreover, Faber showed \cite[Theorem 3.4]{Fa} that  $\vg \geq c \delta_0 (\Gamma)+ \frac{4}{3} \delta_1 (\Gamma) $ with $c=\frac{2}{81}$. Later, the author in \cite[Theorem 2.11]{C5} improved this result by showing that $c$ can be taken as $ \frac{1}{39}$. However, these lower bounds are not sharp.

Contributions of this paper are as follows:

(i) We determine all polarized metrized graphs $(\ga,\bq)$ of genus $3$.

(ii) We explicitly compute various invariants including $\vg$, $\ed$ and $\lag$ of all $(\ga,\bq)$ of genus $3$.

(iii) We give sharp lower bounds for each of $\vg$, $\ed$ and $\lag$. Namely, we show that if $(\ga,\bq)$ is of genus $3$ and of total length $\elg$, then
$\vg \geq \frac{17}{288}\elg$, $\lag \geq \frac{3}{28} \elg$ and $\ed \geq \frac{2}{9} \elg$. As a result, via Zhang's work, this improves the bound given for Effective Bogomolov Conjecture for genus $3$ case (see \cite[Theorems 2.3 and 2.4]{C5}). Note that the bound $\frac{17}{288}\elg$ for $\vg$ was conjectured by X. Faber in \cite[Remark 5.1]{Fa}. This lower bound is attained when $\ga$ is a regular tetrahedral graph. We note that the number $288$ also appears in Tartaglia's formula for the volume of a tetrahedron,
which is a generalization of Heron's formula for the area of a triangle.

(iv) We obtain, as a by product, a highly nontrivial inequality that holds for any nonnegative six real numbers (see inequalities (\ref{eqn main inequality}), (\ref{eqn main inequality2}) and (\ref{eqn main inequality0}) below). Interestingly, the terms that appear in these inequality corresponds to certain cycles of a tetrahedral graph which is not necessarily in $\RR^3$.

\vskip 0.5 in
\section{pm-graphs and their invariants}\label{sec pm-graphs}

In this section, we first give brief descriptions of a metrized graph $\ga$, a polarized metrized graph $(\ga,\bq)$, invariants $\tg$, $\tcg$, $\ed$, $\vg$, $\lag$ and $Z(\ga)$ associated to $(\ga,\bq)$.

A metrized graph $\ga$ is a finite connected graph equipped with a distinguished parametrization of each of its edges.
A metrized graph $\ga$ can have multiple edges and self-loops.
For any given $p \in \ga$,
the number $\va(p)$ of directions emanating from $p$ will be called the \textit{valence} of $p$.
By definition, there can be only finitely many $p \in \ga$ with $\va(p)\not=2$.

For a metrized graph $\ga$, we will denote a vertex set for $\ga$ by $\vv{\ga}$.
We require that $\vv{\ga}$ be finite and non-empty and that $p \in \vv{\ga}$ for each $p \in \ga$ if $\va(p)\not=2$. For a given metrized graph $\ga$, it is possible to enlarge the
vertex set $\vv{\ga}$ by considering additional valence $2$ points as vertices.

For a given metrized graph $\ga$ with vertex set $\vv{\ga}$, the set of edges of $\ga$ is the set of closed line segments with end points in $\vv{\ga}$. We will denote the set of edges of $\ga$ by $\ee{\ga}$. However, if
$e_i$ is an edge, by $\ga-e_i$ we mean the graph obtained by deleting the interior of $e_i$.

We define the genus of $\ga$ to be the first Betti number $g(\ga):=e-v+1$ of the graph $\ga$, where $e$ and $v$ are the number of edges and vertices of $\ga$, respectively.

Length of an edge of $\ga$ is a positive real number. If we denote the length of an edge $e_i \in \ee{\ga}$ by $\li$, the total length of $\ga$, which is denoted by $\elg$, is given by $\elg=\sum_{i=1}^e\li$.


The tau constant $\tg$ of a metrized graph $\ga$ was initially defined by Baker and Rumely in \cite[Section 14]{BRh}.
The following lemma gives a description of the tau constant. In particular, it implies that the tau constant is positive.
\begin{lemma}\cite[Lemma 14.4]{BRh}\label{lemtauformula}
Let $r(x,y)$ be the resistance function on $\ga$. For any fixed $y$ in $\ga$,
$\tg =\frac{1}{4}\int_{\ga}\big(\frac{\partial}{\partial x} r(x,y) \big)^2dx$.
\end{lemma}
One can find more detailed information on $\tg$ in articles \cite{C1}, \cite{C2}, \cite{C3} and \cite{C6}.
For more information about the resistance function $r(x,y)$ on a metrized graph, one can consult to
the articles \cite{BRh}, \cite{BF} and \cite{C2}.


Let $\ga$ be a metrized graph with a vertex set $\vv{\ga}$, and let $\bq : \ga \rightarrow  \NN$
be a function supported on a subset of $\vv{\ga}$. That is, $\bq(s)=0$ for all $s \in \ga-\vv{\ga}$, and whenever $\bq(s) > 0$ we must have $s \in \vv{\ga}$.

A divisor on $\ga$ is a formal sum $\sum n_i p_i$, where $a_i \in \ZZ$ and $p_i \in \ga$ for every $i$.
A divisor $\sum n_i p_i$ on $\ga$ is called effective if $n_i \geq 0$ for all $i$.

The canonical divisor $K$ of $(\ga,\bq)$ is defined as follows:
\begin{equation}\label{eqn app2a}
\begin{split}
K  = \sum_{p \in \vv{\ga}} (\va(p)-2+2 \bq (p))p.
\end{split}
\end{equation}
The pair $(\ga,\bq)$ is called a \textit{polarized metrized graph} (pm-graph in short) if $K$ is an effective divisor. Whenever $\bq=0$, $(\ga,\bq)$ is called a simple pm-graph. We define the \textit{genus} $\gc (\ga)$ of a pm-graph $(\ga,\bq)$ as follows:
\begin{equation}\label{eqn genus}
\begin{split}
\gc (\ga) =g(\ga)+\sum_{p \in \vv{\ga}}\bq (p).
\end{split}
\end{equation}
If $\ga$ under consideration is clear, we simply use notations $g$ and $\gc$ instead of $\gc(\ga)$ and $g(\ga)$, respectively.
\begin{remark}\label{rem effective}
For each $p \in \vv{\ga}$, $\va(p)-2+2 \bq (p) \geq 0$ and $\bq (p)\geq 0$, since the canonical divisor
$K$ is effective and $\bq$ is nonnegative. In particular, if $\va(p)=1$ for some $p \in \ga$, we should have $p \in \vv{\ga}$ and $\bq (p) \geq 1$.
\end{remark}


On a pm-graph $(\ga,\bq)$, we defined and studied the invariant $\tcg$ in \cite{C1} and \cite{C5} as follows:
\begin{equation}\label{eqn tcg}
\begin{split}
\tcg=\sum_{p, \, q \in \, \vv{\ga}}(\va(p)-2+2 \bq (p))(\va(q)-2+2 \bq (q))r(p,q).
\end{split}
\end{equation}
We have $\tcg \geq 0$ for any pm-graph $\ga$, since the canonical divisor $K$ is effective.

Let $\nn{ad}(x)$ be the admissible measure associated to $K$ (defined by Zhang \cite[Lemma 3.7]{Zh1}).
Next, we give definitions of the invariants $\ed$, $\vg$, $\lag$ and $Z(\ga)$ (c.f. \cite[Section 4.1]{Zh2}) of $\ga$:
\begin{equation}\label{eqn app3}
\begin{aligned}
\ed &=\iint_{\Gamma \times \Gamma} r(x,y) \dd{K}(x)  \nn{ad}(y), & \quad Z(\ga) &=\frac{1}{2} \iint_{\Gamma \times \Gamma} r(x,y) \nn{ad}(x) \nn{ad}(y),
\\
\vg &=3 \gc \cdot Z(\ga) -\frac{1}{4} (\ed +\ell (\ga)), & \quad \lag &=\frac{\gc-1}{6 (2 \gc+1)} \vg +\frac{1}{12}(\ed +\elg).
\end{aligned}
\end{equation}

We can express each invariant given in \eqnref{eqn app3} in terms of $\tg$ and $\tcg$ (\cite[Propositions 4.6, 4.7, 4.9 and Theorem 4.8]{C5}):
\begin{theorem}\label{thm pmginv and tau}
Let $(\ga,\bq)$ be a pm-graph with $\gc =3$. Then we have
\begin{align*}\label{eqn pmginv and tau}
\vg &= \frac{13}{3}\tg +\frac{\tcg}{12}-\frac{\elg}{4},&  Z(\ga) &= \frac{5}{9} \tg+\frac{\tc(\ga)}{72},
\\ \lag &= \frac{3}{7} \tg+\frac{\tcg}{56}+\frac{\elg}{14},&  \ed & = \frac{8}{3} \tg + \frac{\tcg}{6}.
\end{align*}
%
%
\end{theorem}

Let $p$ be a point in a pm-graph $\ga$ such that $p \not \in \vv{\ga}$. That is, $p$ is not an end point of any edge in $\ee{\ga}$. Let $\ga$-${p}$ be the pm-graph obtained from $\ga$ by removing $p$ and adding two points $p_1$ and $p_2$ to make the remaining parts closed. That is, $\ga$ is obtained
from $\ga$-$p$ by identifying the end points $p_1$ and $p_2$. Following Zhang's definition \cite[Section 4.1]{Zh2}, we call $p$ is of type $0$ if $\ga$-${p}$ is connected. This happens when $p$ is contained in an edge such that removing the edge does not disconnect $\ga$. If $p$ is not of type $0$, $\ga$-${p}$ is a union of two connected metrized subgraphs with functions $\bq_1$ and $\bq_2$ that are restrictions of $\bq$ and satisfy $\bq_1(p_1)=\bq_2(p_2)=0$. By applying \eqnref{eqn genus}, we see that the subgraphs are of genus $i$ and $\gc - i$ for some integer $i \in (0,\gc/2]$. In this case, we call $p$ is of type $i$. For each integer $i \in [0,\gc/2]$, let $\ga_i$ be the subgraph of $\ga$ of points of type $i$, and let $\ell_i(\ga)$ be the total length of $\ga_i$. We use the invariants $\delta_i (\Gamma):=\ell_i(\ga)$ for each $i \geq 0$ (see \cite{C5} or \cite{Zh2} for geometric meaning of these invariants). Therefore, whenever $\gc =3$, we have only type $0$ and $1$ points, so we consider only the invariants $\delta_0 (\Gamma)$ and $\delta_1 (\Gamma)$ for which we have $\elg=\delta_0 (\Gamma)+\delta_1 (\Gamma)$.

\begin{remark}\label{rem valence}
Given a pm-graph $(\ga,\bq)$ with a vertex set $\vv{\ga}$ containing at least two elements, suppose $\bq (s)=0$ and $\va(s)=2$ for some $s \in \vv{\ga}$, then removing $s$ from the vertex set of $\ga$ does not change $\tcg$. Similarly, if a vertex $s$ is such that $\va(s)=2$ and $\vv{\ga}-\{s \}$ has at least one element, then removing $s$ from $\vv{\ga}$ does not change $\tg$.
(such vertices are called eliminable vertices in \cite[pg. 152]{KY1}). We call these the {\em valence property} of $\tg$ and $\tcg$, see \cite[Remark 2.10]{C2}).
Therefore, $\ed$, $Z(\ga)$, $\vg$ and $\lag$ do not change under this process by \thmref{thm pmginv and tau}. That is, each of these invariants has the valence property \cite[Remark 2.4]{ZCGraphInvComp}.
\end{remark}

\remref{rem valence} is very helpful to determine the possible pm-graphs of a given genus.
As long as $\vv{\ga}$ is non-empty, it will be enough to consider pm-graphs not having any vertex $p$ with $\va(p)=2$ and $\bq(p)=0$.
In fact, in this way we choose only one model of a given pm-graph among all the equivalent models.

Recall that we consider only pm-graphs with $\gc =3$ in this paper. Such pm-graphs can have any $g \in \{0,\, 1, \, 2, \, 3 \}$ by \eqnref{eqn genus}. We designate a section for each such value of $g$ in the rest of the article. In each case, we first determine the pm-graphs with the desired $g$ and $\gc$ up to equivalence. Then, we compute all the relevant invariants. Finally, we find sharp lower bounds to the invariants $\tg$ and $\vg$, $\lag$ for all $\ga$ under consideration. We can use the techniques developed in \cite{C2} and \cite{C1} along with \thmref{thm pmginv and tau} to compute these invariants. This is what we did in this paper. Alternatively, we can compute these invariants by using the algorithms given in \cite{C6} and \cite{ZCGraphInvComp}. For example, how we compute the tau constant for the metrized graph in part $IX$ of \figref{fig caseIV} is illustrated in \cite[Example 5.2]{C6}, and computation of invariants of the pm-graph
in part $XIV$ of \figref{fig caseIV} are done in \cite[Example 1]{ZCGraphInvComp}.

Note that we could consider only pm-graphs $XIII$ and $XIV$ in \figref{fig caseIV} to obtain lower bound of $\vg$ for all pm-graphs of $\gc=3$. The proof of this fact was given in \cite[pages 360 and 365]{Fa} (see also \cite[pages 549 and 550]{C5}). Another approach utilizing this fact was also known to Yamaki (see \cite[page 67]{KY3} and \cite[page 160]{KY1}). However, our aim in this paper is more than finding the sharp lower bound, so we worked on all possible pm-graphs.

\section{The case $g(\ga) = 0$}\label{sec caseI}

Suppose $\ga$ is a pm-graph with $g(\ga)=0$ and $\gc (\ga)=3$. Then \eqnref{eqn genus} becomes $3=\sum_{p \in \vv{\ga}}\bq (p)$.
Moreover, it follows from \remref{rem effective} that such a pm-graph can have at most three vertices with valence exactly $1$. On the other hand, $g(\ga)=0$ implies $e=v-1$.

\begin{figure}
\centering
\includegraphics[scale=0.65]{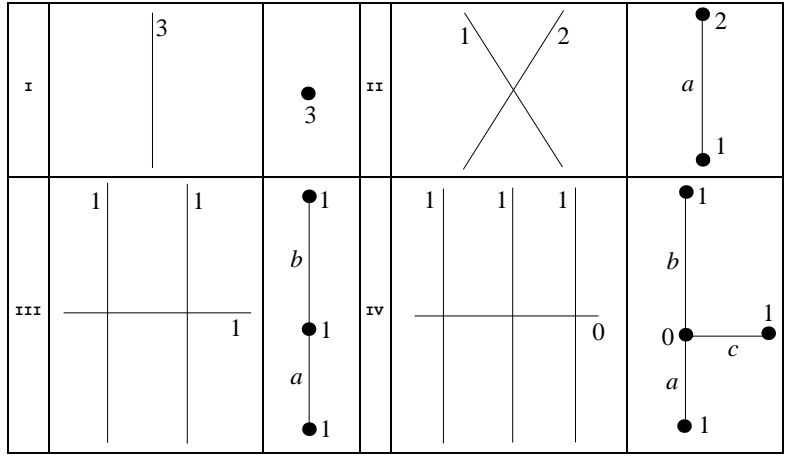}  \caption{Irreducible components of genus 3 curves and their dual graphs when the dual graphs have genus $0$, i.e. when $\gc=3$ and $g=0$.} \label{fig caseI}
\end{figure}

If $\ga$ has no edges, then it is just a point $p$ with $\bq (p)=3$, as shown in part $I$ of \figref{fig caseI}. Otherwise, pm-graphs $(\ga,\bq)$ in this case are tree pm-graphs, which have at least two vertices with $\va(p)=1$. \figref{fig caseI} illustrates the possible pm-graphs satisfying these conditions.

We have $\tg=\frac{\elg}{4}$ since $\ga$ is a metrized graph that is a tree graph \cite[Equation 14.3]{BRh}, and we use \eqnref{eqn tcg} to compute $\tcg$. Then we compute $\vg$, $\lag$ and $\ed$ by using \thmref{thm pmginv and tau}. The results are given in \tabref{tab caseI}.

\begin{table}
\begin{center}
\begin{tabular}{|c|c|c|c|c|c|c|c|c|}
  \hline
 $$ \T \B & $\elg$ & $\delta_0(\ga)$& $\delta_1(\ga)$ & $\tg$ & $\tcg$ & $\vg$ &  $\lag$ & $\ed$ \\
 \hline
 $I$ \T \B & $0$& $0$& $0$ &$0$ & $0$ & $0$ & $0$ & $0$\\
 \hline
 $II$ \T \B & $a$& $0$& $\elg$ &$\frac{\elg}{4}$ & $6\elg$ & $\frac{4\elg}{3}$ & $\frac{2\elg}{7}$ & $\frac{5\elg}{3}$\\
 \hline
 $III$ \T \B & $a+b$ & $0$& $\elg$ &$\frac{\elg}{4}$ & $6\elg$ & $\frac{4\elg}{3}$ & $\frac{2\elg}{7}$& $\frac{5\elg}{3}$\\
 \hline
 $IV$ \T \B & $a+b+c$ & $0$& $\elg$ &$\frac{\elg}{4}$ & $6\elg$ & $\frac{4\elg}{3}$ & $\frac{2\elg}{7}$ & $\frac{5\elg}{3}$\\
 \hline
\end{tabular}
\end{center}  \caption{Pm-graph invariants when $g(\ga)=0$} \label{tab caseI}
\end{table}
We exclude the case $I$ in \tabref{tab caseI} as $\elg =0$. In the other three cases, we have  $\vg = \frac{4}{3} \elg$, $\lag=\frac{2\elg}{7}$ and $\ed=\frac{5\elg}{3}$.

\section{The case $g(\ga) = 1$}\label{sec caseII}

\begin{figure}
\centering
\includegraphics[scale=0.62]{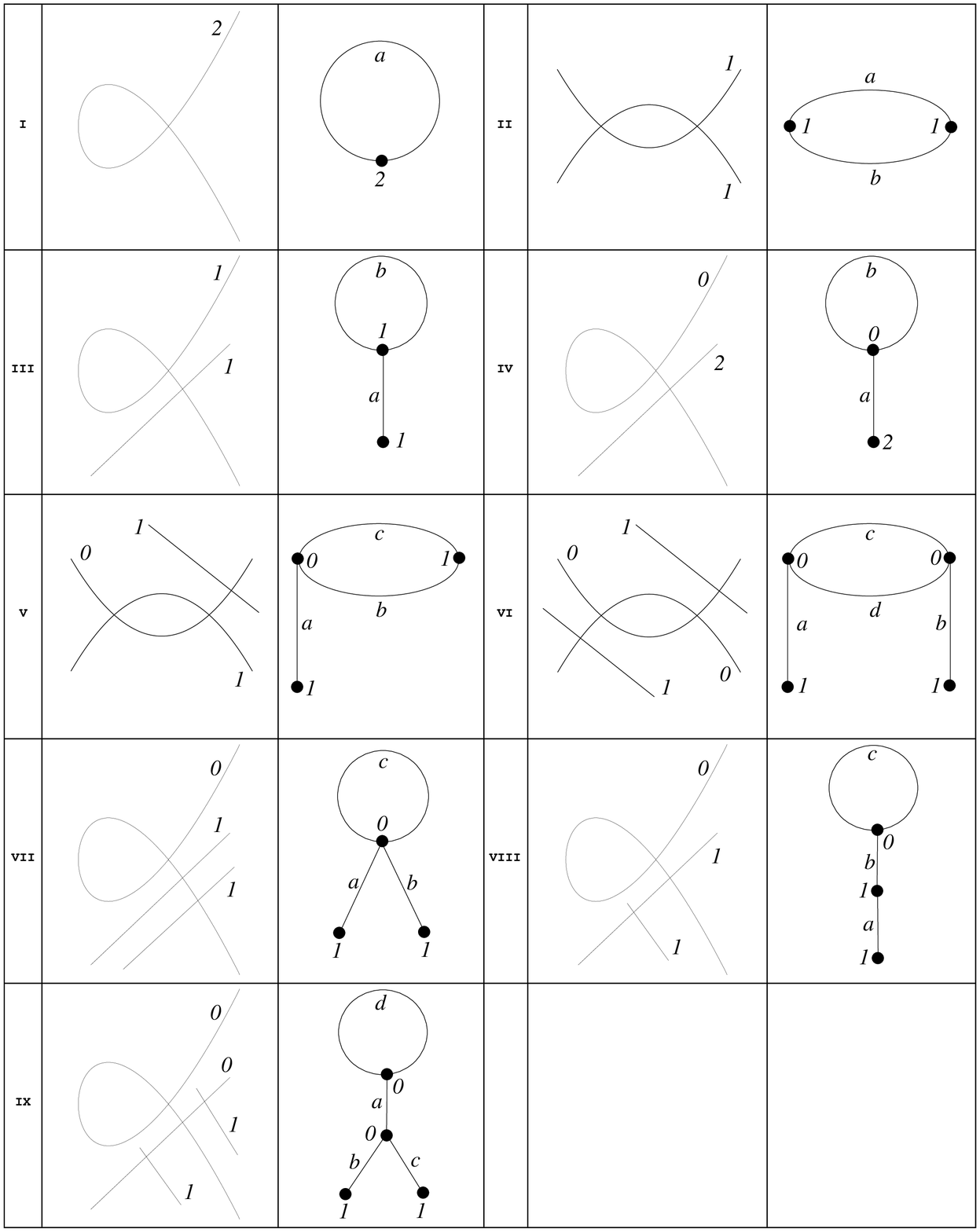} \caption{Irreducible components of genus 3 curves and their dual graphs when the dual graphs have genus $1$, i.e. when $\gc=3$ and $g=1$.} \label{fig caseII}
\end{figure}

In this section, we consider pm-graphs $(\ga,\bq)$ with $g=1$ and $\gc =3$. It follows from \eqnref{eqn genus}
that $2=\sum_{p \in \vv{\ga}}\bq (p)$. Thus, such a pm-graph can have at most two vertices with valence exactly $1$ by \remref{rem effective}. Since $g=1$, we have $e=v$. Based on these observations,
\figref{fig caseII} illustrates the possible pm-graphs satisfying these conditions.

We have $\tg=\frac{\elg}{4}$ when $\ga$ is a tree metrized graph \cite[Equation 14.3]{BRh}. Moreover, $\tg=\frac{\elg}{12}$ when $\ga$ is a circle metrized graph \cite[Corollary 2.17]{C2}. We use these facts and additive property of tau constant \cite[page 15]{C2} to compute $\tg$ for each pm-graphs listed in \figref{fig caseII}. We again use \eqnref{eqn tcg} to compute $\tcg$. Then we compute $\vg$, $\lag$ and $\ed$ by using \thmref{thm pmginv and tau}.
The invariants $\delta_0(\ga)$ and $\delta_1(\ga)$ are determined by using their definitions and by considering the topology of $\ga$.
The results are given in \tabref{tab caseIIa} and \tabref{tab caseIIb}. As can be seen from \tabref{tab caseIIb}, we have $\vg \geq \frac{1}{9} \elg$, $\lag \geq \frac{3\elg}{28}$, $\ed \geq \frac{2\elg}{9}$, and these lower bounds are attained by the pm-graph given in part $I$ of \figref{fig caseII}.

\begin{table}
\begin{center}
\begin{tabular}{|c|c|c|c|c|c|}
  \hline
 $$ \T \B & $\elg$ & $\delta_0(\ga)$& $\delta_1(\ga)$ & $\tg$ & $\tcg$ \\
 \hline
 $I$ \T \B & $a$& $a$ & $0$ & $\frac{\elg}{12}$ & $0$\\
 \hline
 $II$ \T \B & $a+b$& $a+b$ & $0$ & $\frac{\elg}{12}$ & $\frac{8ab}{a+b}$\\
 \hline
 $III$ \T \B & $a+b$ & $b$ & $a$ & $\frac{\elg}{12}+\frac{a}{6}$ & $6 a$\\
 \hline
 $IV$ \T \B & $a+b$ & $b$ & $a$ & $\frac{\elg}{12}+\frac{a}{6}$ & $6 a$ \\
 \hline
$V$ \T \B & $a+b+c$ & $b+c$ & $a$ & $\frac{\elg}{12}+\frac{a}{6}$ & $6 a+\frac{8 b c}{b+c}$ \\
 \hline
 $VI$ \T \B & $a+b+c+d$ & $d+c$ & $a+b$ & $\frac{\elg}{12}+\frac{a+b}{6}$ & $6(a+b)+\frac{8 c d}{c+d}$ \\
 \hline
$VII$ \T \B & $a+b+c$ & $c$ & $a+b$ &$\frac{\elg}{12}+\frac{a+b}{6}$ & $6 (a+b)$ \\
 \hline
$VIII$ \T \B & $a+b+c$ & $c$ & $a+b$ &$\frac{\elg}{12}+\frac{a+b}{6}$ & $6 (a+b)$ \\
 \hline
 $IX$ \T \B & $a+b+c+d$ & $d$ & $a+b+c$ &$\frac{\elg}{12}+\frac{a+b+c}{6}$ & $6 (a+b+c)$ \\
 \hline
\end{tabular}
\end{center}  \caption{Pm-graph invariants when $g(\ga)=1$, part 1.} \label{tab caseIIa}
\end{table}

\begin{table}
\begin{center}
\begin{tabular}{|c|c|c|c|c|c|c|}
  \hline
 $$ \T \B &  $\vg$ &  $\lag$ & $\ed$ \\
 \hline
 $I$ \T \B &  $\frac{\elg}{9}$ & $\frac{3\elg}{28}$ & $\frac{2\elg}{9}$ \\
 \hline
 $II$ \T \B &  $\frac{\elg}{9}+\frac{2ab}{3 (a + b)}$ & $\frac{3\elg}{28}+\frac{ab}{7 (a + b)}$ & $\frac{2\elg}{9}+\frac{4 a b}{3(a+b)}$ \\
 \hline
 $III$ \T \B &  $\frac{\elg}{9}+\frac{11a}{9}$ & $\frac{3\elg}{28}+\frac{5a}{28}$& $\frac{2\elg}{9}+\frac{13a}{9}$ \\
 \hline
 $IV$ \T \B &  $\frac{\elg}{9}+\frac{11a}{9}$ & $\frac{3\elg}{28}+\frac{5a}{28}$& $\frac{2\elg}{9}+\frac{13a}{9}$ \\
 \hline
$V$ \T \B &  $\frac{\elg}{9}+\frac{6 b c + 11 a (b + c)}{9 (b + c)}$ & $\frac{3\elg}{28}+\frac{4b c +5 a (b + c)}{28 (b + c)}$ & $\frac{2\elg}{9}+\frac{12 b c+ 13 a (b + c)}{9 (b + c)}$ \\
 \hline
 $VI$ \T \B &  $\frac{\elg}{9}+\frac{6 c d + 11( a+ b) (c + d)}{9 (c + d)}$ & $\frac{3\elg}{28}+\frac{4cd+5(a +b)(c + d) }{28 (c + d)}$ & $\frac{2\elg}{9}+\frac{12cd+13 (a +  b)(c+d)}{9(c + d)}$ \\
 \hline
$VII$ \T \B &  $\frac{\elg}{9}+\frac{11 (a + b)}{9}$ & $\frac{3\elg}{28}+\frac{5 (a + b)}{28}$ & $\frac{2\elg}{9}+\frac{13 (a + b))}{9}$  \\
 \hline
$VIII$ \T \B &  $\frac{\elg}{9}+\frac{11 (a + b)}{9}$ & $\frac{3\elg}{28}+\frac{5 (a + b)}{28}$ & $\frac{2\elg}{9}+\frac{13 (a + b)}{9}$  \\
 \hline
 $IX$ \T \B &  $\frac{\elg}{9}+\frac{11 (a + b+c)}{9}$ & $\frac{3\elg}{28}+\frac{5 (a + b+c)}{28}$ & $\frac{2\elg}{9}+\frac{13 (a + b+c)}{9}$  \\
 \hline
\end{tabular}
\end{center}  \caption{Pm-graph invariants when $g(\ga)=1$, part 2.} \label{tab caseIIb}
\end{table}

\section{The case $g(\ga) = 2$}\label{sec caseIII}

\begin{figure}
\centering
\includegraphics[scale=0.59]{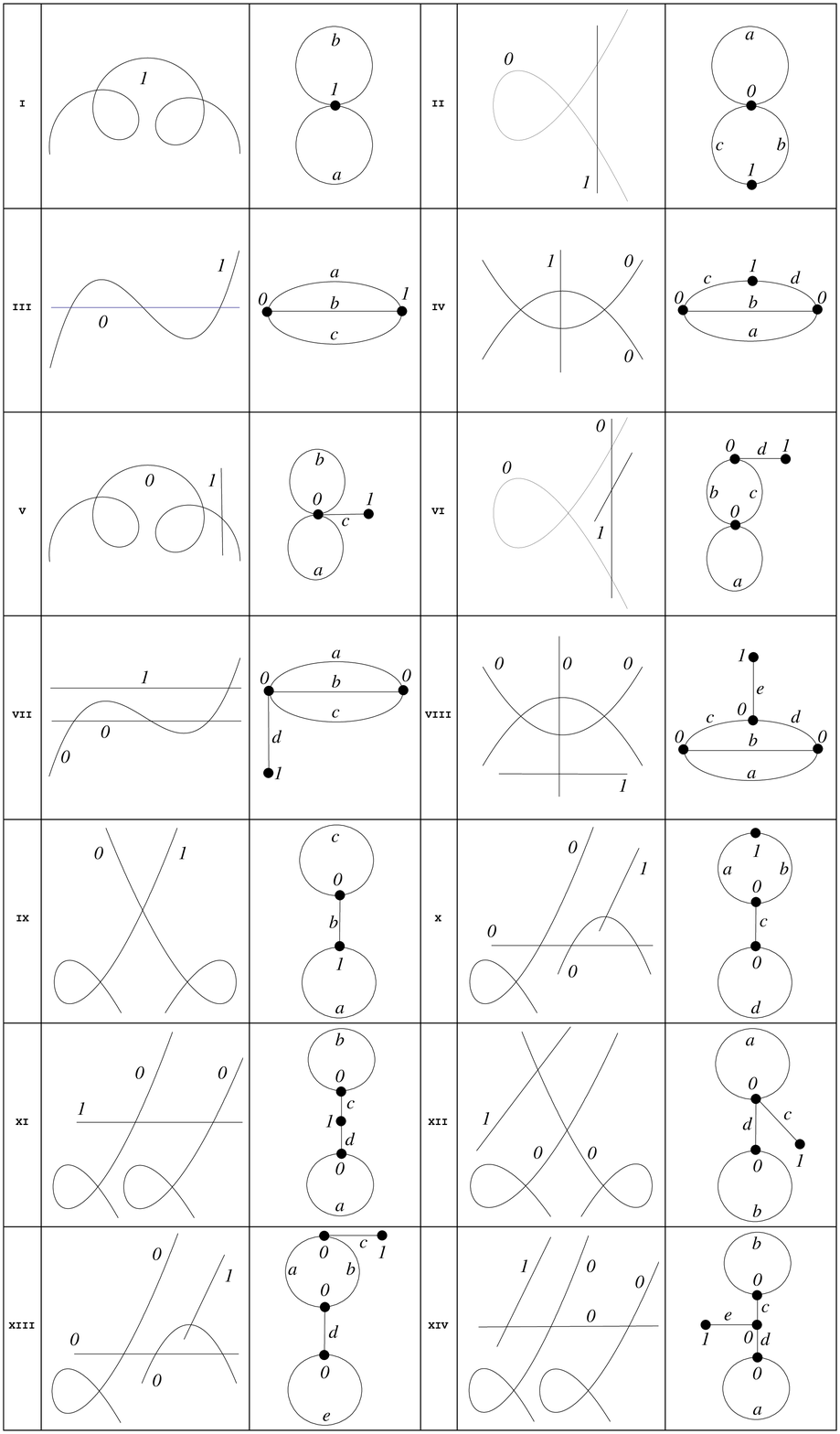} \caption{Irreducible components of genus 3 curves and their dual graphs when the dual graphs have genus $2$, i.e. when $\gc=3$ and $g=2$.} \label{fig caseIII}
\end{figure}

In this section, we consider pm-graphs $(\ga,\bq)$ with $g=2$ and $\gc =3$. Using \eqnref{eqn genus} we see that
$\bq (p)=1$ for only one vertex $p \in \vv{\ga}$ and that $\bq (p)=0$ for all the remaining vertices. By \remref{rem effective} again, such a pm-graph can have at most one vertex with valence exactly $1$. Moreover, we have $e=v+1$. Based on these observations,
\figref{fig caseIII} illustrates all the pm-graphs satisfying these conditions.

A metrized graph with two vertices and $m$ multiple edges connecting these two vertices is called $m$-banana. For such $\ga$ we know how to compute $\tg$ \cite[Proposition 8.3]{C2}. Then using tau formulas for tree and circle metrized graphs along with the additive property, one can compute $\tg$ for each of the pm-graphs given in \figref{fig caseIII}.
Again we use its definition to compute $\tcg$. Then we compute the remaining invariants $\vg$, $\lag$ and $\ed$ by using \thmref{thm pmginv and tau}.
The invariants $\delta_0(\ga)$ and $\delta_1(\ga)$ are determined by using their definitions and by considering the topology of $\ga$.
The results are given in \tabref{tab caseIIIa} and \tabref{tab caseIIIb}.

As can be seen from the values of $\vg$, we have $\vg \geq \frac{1}{9} \elg$ if $\ga$ is one of the pm-graphs given in parts $I$, $II$,
$V$, $VI$, $IX$, $X$, $XI$, $XII$, $XIII$, $XIV$.

For the pm-graph $\ga$ of type $III$, we see that $\frac{a+b+c}{3}-\frac{3}{1/a+1/b+1/c} \geq 0$ by Arithmetic-Harmonic Mean inequality. Note that $\elg = a+b+c$. Therefore, we have $\vg \geq \frac{1}{9} \elg - \frac{2}{9} \frac{1}{9} \elg = \frac{7}{81} \elg$, which is the sharp lower bound for this type of pm-graphs, because  $\vg = \frac{7}{81} \elg$ whenever $a=b=c$.

Using the same inequality $\frac{a+b+c}{3}-\frac{3}{1/a+1/b+1/c} \geq 0$, we see that  $\vg \geq \frac{7}{81} \elg$ for the pm-graph of type $VII$.

Similarly, if we use Arithmetic-Harmonic Mean inequality for $a$, $b$ and $c+d$, we again obtain that $\vg \geq \frac{7}{81} \elg$ for the pm-graphs as in type $IV$ and $VIII$.

In any case,  we have the sharp lower bound $\vg \geq \frac{7}{81} \elg$ whenever $g = 2$.

Using the results from \tabref{tab caseIIIb} we have, as in the previous section, $\lag \geq \frac{3\elg}{28}$ and $\ed \geq \frac{2\elg}{9}$, and these lower bounds are attained by the pm-graph given in part $I$ of \figref{fig caseIII}.

\begin{table}
\begin{center}
\begin{tabular}{|c|c|c|c|c|}
  \hline
 $$ \T \B & $\elg$ & $\delta_1(\ga)$ & $\tg$ & $\tcg$ \\
 \hline
 $I$ \T \B & \small{$a+b$}  & $0$ & $\frac{\elg}{12}$ & $0$\\
 \hline
 $II$ \T \B & \small{$a+b+c$}  & $0$ & $\frac{\elg}{12}$ & $\frac{8bc}{b+c}$\\
 \hline
 $III$ \T \B & \small{$a+b+c$}  & $0$ & $\frac{\elg}{12}-\frac{abc}{6(ab +ac+bc)}$ & $\frac{6abc}{ab +ac+bc}$\\
 \hline
 $IV$ \T \B & \small{$a+b+c+d$}  & $0$ & $\frac{\elg}{12}-\frac{ab(c+d)}{6(ab +(a+b)(c+d))}$ & $\frac{6ab(c+d)+8(a+b)cd}{ab +(a+b)(c+d)}$ \\
 \hline
$V$ \T \B & \small{$a+b+c$}  & $c$ & $\frac{\elg}{12}+\frac{c}{6}$ & $6c$ \\
 \hline
 $VI$ \T \B & \small{$a+b+c+d$}  & $d$ & $\frac{\elg}{12}+\frac{d}{6}$ & $6d+\frac{8 b c}{b+c}$ \\
 \hline
$VII$ \T \B & \small{$a+b+c+d$}  & $d$ &$\frac{\elg}{12}+\frac{d}{6}-\frac{abc}{6(ab +ac+bc)}$ & $6d+\frac{6abc}{ab +ac+bc}$ \\
 \hline
\small{$VIII$} \T \B & \small{$a+b+c+d+e$}  & $e$ &$\frac{\elg}{12}+\frac{e}{6}-\frac{ab(c+d)}{6(ab +(a+b)(c+d))}$ & $6e+\frac{6ab(c+d)+8(a+b)cd}{ab +(a+b)(c+d)}$ \\
 \hline
 $IX$ \T \B & \small{$a+b+c$} & $b$ &$\frac{\elg}{12}+\frac{b}{6}$ & $6 b$ \\
 \hline
 $X$ \T \B & \small{$a+b+c+d$}  & $c$ &$\frac{\elg}{12}+\frac{c}{6}$ & $6c+\frac{8ab}{a+b}$ \\
 \hline
 $XI$ \T \B & \small{$a+b+c+d$}  & $c+d$ &$\frac{\elg}{12}+\frac{c+d}{6}$ & $6(c+d)$ \\
 \hline
 \small{$XII$} \T \B & \small{$a+b+c+d$}  & $c+d$ &$\frac{\elg}{12}+\frac{c+d}{6}$ & $6(c+d)$ \\
 \hline
 \small{$XIII$} \T \B & \small{$a+b+c+d+e$}  & $d+c$ &$\frac{\elg}{12}+\frac{c+d}{6}$ & $6(c+d)+\frac{8ab}{a+b}$ \\
 \hline
 \small{$XIV$} \T \B & \small{$a+b+c+d+e$}  & $c+d+e$ &$\frac{\elg}{12}+\frac{c+d+e}{6}$ & $6(c+d+e)$ \\
 \hline
\end{tabular}
\end{center}  \caption{Pm-graph invariants when $g(\ga)=2$, part 1. We have $\delta_0(\ga)=\elg-\delta_1(\ga)$, and $\delta_i(\ga)=0$ for all $i \geq 2$.} \label{tab caseIIIa}
\end{table}

\begin{table}
\begin{center}
\begin{tabular}{|c|c|c|c|c|c|c|}
  \hline
 $$ \T \B &  $\vg$ &  $\lag$ & $\ed$ \\
 \hline
 $I$ \T \B &  $\frac{\elg}{9}$ & $\frac{3\elg}{28}$ & $\frac{2\elg}{9}$ \\
 \hline
 $II$ \T \B &  $\frac{\elg}{9}+\frac{2b c}{3(b+c)}$ & $\frac{3\elg}{28}+\frac{b c}{7(b+c)}$ & $\frac{2\elg}{9}+\frac{4b c}{3(b+c)}$ \\
 \hline
 $III$ \T \B &  $\frac{\elg}{9}-\frac{2a b c}{9(ab+ac+bc)}$ & $\frac{3\elg}{28}+\frac{a b c}{28(ab+ac+bc)}$& $\frac{2\elg}{9}+\frac{5a b c}{9(ab+ac+bc)}$ \\
 \hline
 $IV$ \T \B &  $\frac{\elg}{9}+\frac{6cd(a+b)-2ab(c+d)}{9(ab+(a+b)(c+d))}$ & $\frac{3\elg}{28}+\frac{4cd(a+b)+ab(c+d)}{28(ab+(a+b)(c+d))}$& $\frac{2\elg}{9}+\frac{12cd(a+b)+5ab(c+d)}{9(ab+(a+b)(c+d))}$ \\
 \hline
$V$ \T \B &  $\frac{\elg}{9}+\frac{11c}{9}$ & $\frac{3\elg}{28}+\frac{5c}{28}$ & $\frac{2\elg}{9}+\frac{13c}{9}$ \\
 \hline
 $VI$ \T \B &  $\frac{\elg}{9}+\frac{6bc+11d(b+c)}{9(b+c)}$ & $\frac{3\elg}{28}+\frac{4bc+5d(b+c)}{28(b+c)}$ & $\frac{2\elg}{9}+\frac{12bc+13d(b+c)}{9(b+c)}$ \\
 \hline
$VII$ \T \B &  $\frac{\elg}{9}+\frac{11d}{9}-\frac{2a b c}{9(ab+ac+bc)}$ & $\frac{3\elg}{28}+\frac{5d}{28}+\frac{a b c}{28(ab+ac+bc)}$ & $\frac{2\elg}{9}+\frac{13d}{9}+\frac{5a b c}{9(ab+ac+bc)}$  \\
 \hline
\small{$VIII$} \T \B &
$\frac{\elg+11e}{9}+\frac{6cd(a+b)-2a b (c+d)}{9(ab+(a+b)(c+d))}$ &
 $\frac{3\elg+5e}{28}+\frac{4cd(a+b)+a b (c+d)}{28(ab+(a+b)(c+d))}$ & $\frac{2\elg+13e}{9}+\frac{12cd(a+b)+5a b (c+d)}{9(ab+(a+b)(c+d))}$  \\
 \hline
 $IX$ \T \B &  $\frac{\elg}{9}+\frac{11b}{9}$ & $\frac{3\elg}{28}+\frac{5b}{28}$ & $\frac{2\elg}{9}+\frac{13b}{9}$  \\
 \hline
 $X$ \T \B &  $\frac{\elg}{9}+\frac{6ab+11c(a+b)}{9(a+b)}$ & $\frac{3\elg}{28}+\frac{4ab+5c(a+b)}{28(a+b)}$ & $\frac{2\elg}{9}+\frac{12ab+13c(a+b)}{9(a+b)}$  \\
 \hline
$XI$ \T \B &  $\frac{\elg}{9}+\frac{11(c+d)}{9}$ & $\frac{3\elg}{28}+\frac{5(c+d)}{28}$ & $\frac{2\elg}{9}+\frac{13(c+d)}{9}$  \\
 \hline
$XII$ \T \B &  $\frac{\elg}{9}+\frac{11(c+d)}{9}$ & $\frac{3\elg}{28}+\frac{5(c+d)}{28}$ & $\frac{2\elg}{9}+\frac{13(c+d)}{9}$  \\
 \hline
 \small{$XIII$} \T \B &  $\frac{\elg}{9}+\frac{6ab+11(a+b)(c+d)}{9(a+b)}$ &
 $\frac{3\elg}{28}+\frac{4ab+5(a+b)(c+d)}{28(a+b)}$ & $\frac{2\elg}{9}+\frac{12ab+13(a+b)(c+d)}{9(a+b)}$  \\
 \hline
 \small{$XIV$} \T \B &  $\frac{\elg}{9}+\frac{11(c+d+e)}{9}$ & $\frac{3\elg}{28}+\frac{5(c+d+e)}{28}$ & $\frac{2\elg}{9}+\frac{13(c+d+e)}{9}$  \\
 \hline
\end{tabular}
\end{center}  \caption{Pm-graph invariants when $g(\ga)=2$, part 2.} \label{tab caseIIIb}
\end{table}

\section{The case $g(\ga) = 3$}\label{sec caseIV}

\begin{figure}
\centering
\includegraphics[scale=0.56]{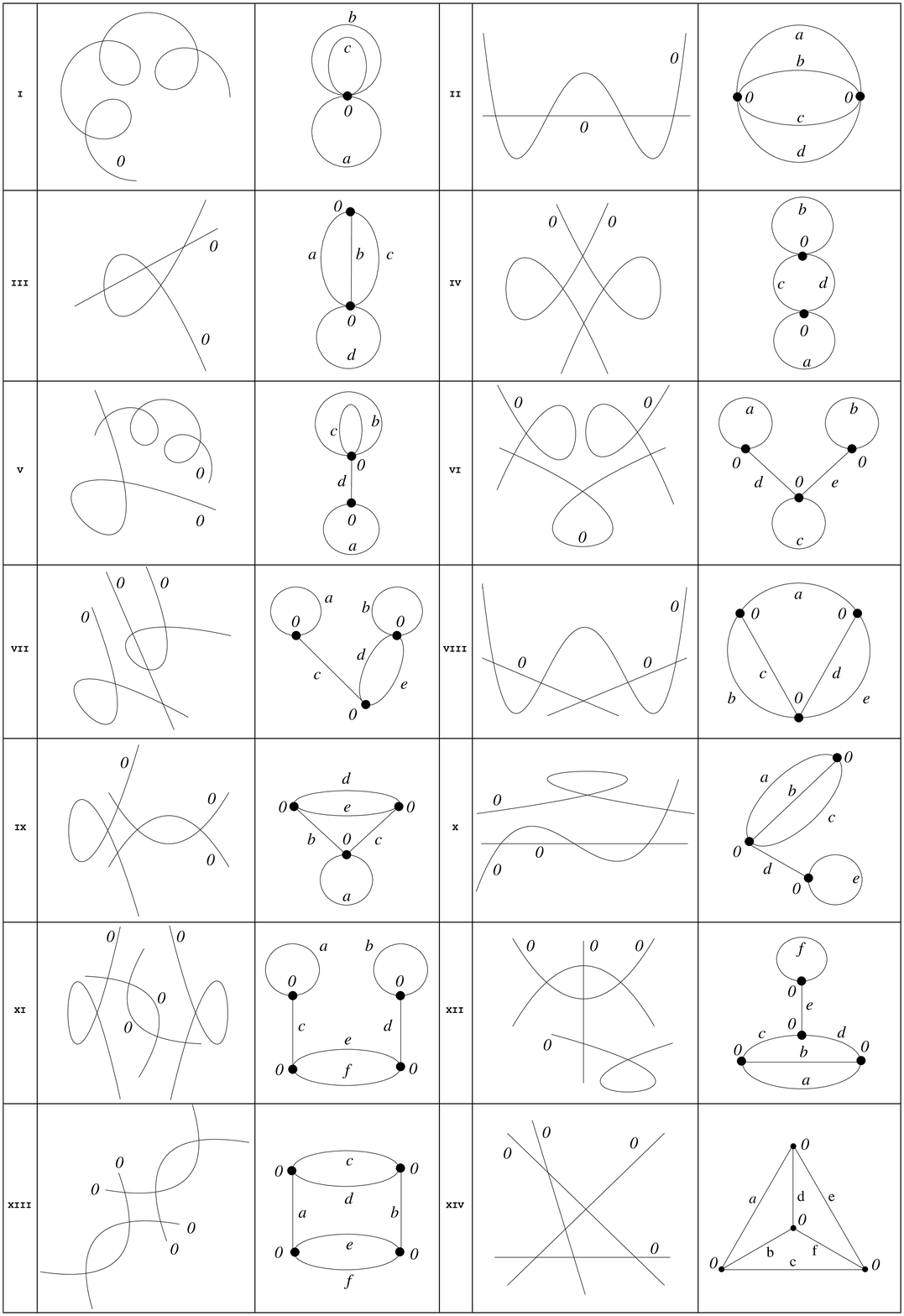} \caption{Irreducible components of genus 3 curves and their dual graphs when the dual graphs have genus 3, i.e. when $\gc=g=3$.} \label{fig caseIV}
\end{figure}

In this section, we consider pm-graphs $(\ga,\bq)$ with $g(\ga) = 3$ and $g=3$. In this case, \eqnref{eqn genus} implies
$\bq (p)=0$ for each vertex $p \in \vv{\ga}$. That is,  $(\ga,\bq)$ is a simple pm-graph.
Using this observation and \remref{rem effective}, we note that $\va(p) \geq 2$ for each $p \in \vv{\ga}$.
Moreover, using \remref{rem valence} we can assume that $\va(p) \geq 3$ for each $p \in \vv{\ga}$ for this section.
By basic graph theory, this implies $e \geq \frac{3}{2}v$. On the other hand, we have $e=v+2$ since $g=3$.
Therefore, we conclude that $1 \leq v \leq 4$ for the simple pm-graphs we can have.
Based on these observations, \figref{fig caseIV} illustrates all the pm-graphs satisfying these conditions.

We compute $\tg$ by using similar techniques as in the previous section except for the pm-graphs in parts $VIII$, $XIII$ and $XIV$ of \figref{fig caseIV}. The simple pm-graph in part $XIV$ is a tetrahedral graph for which we have computed its tau constant in \cite[Example 8.4]{C2}. We can compute the tau constant for the simple pm-graphs in parts  $VIII$ and $XIII$ by using the techniques developed in \cite{C2}, such as \cite[Corollaries 5.3 and 7.4, Propositions 4.6 and 4.5]{C2}.

As in the previous sections, we compute $\tcg$ by using its definition and determining the resistance values between any two vertices of the pm-graph. Note that computing the resistance matrix of the corresponding graph will also help, as it is done in \cite[Example III]{ZCGreenFunction}.

Once the values of $\tg$ and $\tcg$ are obtained, we compute $\vg$, $\lag$ and $\ed$ by using \thmref{thm pmginv and tau}.
As in the previous sections, the topology of $\ga$ is the main factor effecting the invariants $\delta_0(\ga)$ and $\delta_1(\ga)$.

The results for pm-graphs of types $I$-$XII$ in \figref{fig caseIV} are given in \tabref{tab caseIVa}, \tabref{tab caseIVb} and \tabref{tab caseIVc}. Since the values of the invariants are lengthy for the remaining pm-graphs, types $XIII$ and $XIV$, we state them separately in this section.

It is clear from \tabref{tab caseIVc} and the results for pm-graphs of types $XIII$ and $XIV$ that $\lag \geq \frac{3}{28} \elg$ and $\ed \geq \frac{2}{9} \elg$ for all pm-graphs in \figref{fig caseIV}, and this lower bounds are attained for the pm-graph in type $I$.

Clearly, \tabref{tab caseIVb} shows that $\vg \geq \frac{1}{9} \elg$ for pm-graphs of type $I$, $IV$, $V$, $VI$, $VII$, $XI$.

We use Arithmetic-Harmonic Mean inequality for $a$, $b$ ,$c$, i.e.,  $\frac{a+b+c}{3}-\frac{3}{1/a+1/b+1/c} \geq 0$, to derive
$\vg \geq \frac{7}{81} \elg$ for the pm-graphs of types $III$ and $X$. Using the Arithmetic-Harmonic Mean inequality for $b+c$, $d$, $e$ gives  $\vg \geq \frac{7}{81} \elg$ for the pm-graph of type $IX$. Similarly, using the Arithmetic-Harmonic Mean inequality for $a$, $b$, $c+d$ gives  $\vg \geq \frac{7}{81} \elg$ for the pm-graph of type $XII$.

We note that $\vg = \frac{1}{16} \elg + \frac{7}{36}(\frac{a+b+c+d}{4} - \frac{4}{1/a+1/b+1/c+1/d})$ for the pm-graph of type $II$. Therefore, we obtain $\vg \geq \frac{1}{16} \elg$ by using the Arithmetic-Harmonic Mean inequality for $a$, $b$, $c$ and $d$. In particular, $\vg = \frac{1}{16} \elg$ when $a=b=c=d$.

Computations to find the lower bounds of $\vg$ requires more in-depth analysis for pm-graphs of types $VIII$, $XIII$ and $XIV$. Thus, we consider each of these pm-graphs separately.

\textbf{Pm-graphs of type $VIII$:}

Let $H=b^2 c d+b c^2 d+b c d^2+b^2 c e+b c^2 e+b^2 d e-12 b c d e+c^2 d e+b d^2 e+c d^2 e+b c e^2+b d e^2+c d e^2$,
$D = a c e + a b e + c b e + a c d + a b d + c b d + c e d + b e d$, $N=14 a c e + 3 c^2 e + 14 a b e + 3 b^2 e + 3 c e^2 + 3 b e^2 +
  14 a c d + 3 c^2 d + 14 a b d + 3 b^2 d + 3 c d^2 + 3 b d^2$, $M= (c - b)^2 e +  (c - b)^2 d + c (e - d)^2 +
  b (e - d)^2$. Clearly, $D$, $N$ and $M$ are nonnegative.
We note that $H \geq 0$ by applying Arithmetic-Harmonic Mean inequality for $b$, $c$, $d$ and $e$.

Now we note that
$\vg=\frac{1}{16}\elg+\frac{a(N+11 M)+14 H}{288 D}$, which implies that $\vg \geq \frac{1}{16}\elg$. This is the sharp lower bound because
$\vg = \frac{1}{16}\elg$ whenever $a=0$ and $b=c=d=e$.

\begin{table}
\begin{center}
\begin{tabular}{|c|c|c|c|c|}
 \hline
 $$ \T \B & $\elg$ & $\delta_1(\ga)$ & $\tg$  \\
 \hline
 $I$ \T \B & \small{$a+b+c$}  & $0$ & $\frac{\elg}{12}$ \\
 \hline
 $II$ \T \B & \small{$a+b+c+d$}  & $0$ & $\frac{\elg}{12}-\frac{a b c d}{3 (b c d + a (c d + b (c + d)))}$ \\
 \hline
 $III$ \T \B & \small{$a+b+c+d$}  & $0$ & $\frac{\elg}{12}-\frac{abc}{6(ab +ac+bc)}$ \\
 \hline
 $IV$ \T \B & \small{$a+b+c+d$}  & $0$ & $\frac{\elg}{12}$ \\
 \hline
$V$ \T \B & \small{$a+b+c+d$}  & $d$ & $\frac{\elg}{12}+\frac{d}{6}$ \\
 \hline
 $VI$ \T \B & \small{$a+b+c+d+e$}  & $d+e$ & $\frac{\elg}{12}+\frac{d+e}{6}$ \\
 \hline
$VII$ \T \B & \small{$a+b+c+d+e$}  & $c$ &$\frac{\elg}{12}+\frac{c}{6}$ \\
 \hline
\small{$VIII$} \T \B & \small{$a+b+c+d+e$}  & $0$ &$\frac{\elg}{12}-\frac{a b c d + a b c e + a b d e + a c d e + 2 b c d e}{6(
a b d + a c d + b c d + a b e + a c e + b c e + b d e + c d e)}$ \\
 \hline
 $IX$ \T \B & \small{$a+b+c+d+e$} & $0$ &$\frac{\elg}{12}+\frac{b}{6}$ \\
 \hline
 $X$ \T \B & \small{$a+b+c+d+e$}  & $d$ &$\frac{\elg}{12}+\frac{d}{6}-\frac{a b c}{6(a b+a c+b c)}$ \\
 \hline
 $XI$ \T \B & \small{$a+b+c+d+e+f$}  & $c+d$ &$\frac{\elg}{12}+\frac{c+d}{6}$ \\
 \hline
 \small{$XII$} \T \B & \small{$a+b+c+d+e+f$}  & $e$ &$\frac{\elg}{12}+\frac{e}{6}-\frac{a b (c+d)}{6(a b+(a+b) (c+d))}$ \\
 \hline
\end{tabular}
\end{center}  \caption{Values of $\elg$, $\delta_1(\ga)$ and $\tg$ for pm-graphs with $g(\ga)=3$ and $\gc=3$. We have $\delta_0(\ga)=\elg-\delta_1(\ga)$, and $\delta_i(\ga)=0$ for all $i \geq 2$.} \label{tab caseIVa}
\end{table}

\begin{table}
\begin{center}
\begin{tabular}{|c|c|c|}
  \hline
  \T \B & $\tcg$ & $\vg$ \\
 \hline
 $I$ \T \B & $0$  & $\frac{\elg}{9}$ \\
 \hline
 $II$ \T \B & $\frac{8abcd}{b c d + a (c d + b (c + d))}$  & $\frac{\elg}{9}-\frac{7 a b c d}{9 (b c d + a (c d + b (c + d)))}$ \\
 \hline
 $III$ \T \B & $\frac{6abc}{ab +ac+bc}$  & $\frac{\elg}{9}-\frac{2 a b c}{9 (b c + a (b + c))}$ \\
 \hline
 $IV$ \T \B & $\frac{8cd}{c+d}$  & $\frac{\elg}{9}+\frac{2c d}{3(c+d)}$ \\
 \hline
$V$ \T \B &  $6d$ & $\frac{\elg}{9}+\frac{11d}{9}$ \\
 \hline
 $VI$ \T \B & $6(d+e)$  & $\frac{\elg}{9}+\frac{11(d+e)}{9}$ \\
 \hline
$VII$ \T \B & $6c+\frac{8de}{d+e}$  & $\frac{\elg}{9}+\frac{11c}{9}+\frac{2de}{3(d+e)}$ \\
 \hline
\small{$VIII$} \T \B & $\frac{6 a b c d + 6 a b c e + 6 a b d e + 6 a c d e + 8b c d e}{a b d + a c d + b c d + a b e + a c e + b c e + b d e + c d e}$  & $\frac{\elg}{9}-\frac{7 b c d e + 2 a (c d e + b (d e + c (d + e)))}{9 (c d e + a (b + c) (d + e) + b (d e + c (d + e)))}$ \\
 \hline
 $IX$ \T \B & $\frac{8 b c d + 8 b c e + 6 b d e + 6 c d e}{
b d + c d + b e + c e + d e}$  & $\frac{\elg}{9}+\frac{-2 (b + c) d e + 6 b c (d + e)}{9 (d e + (b + c) (d + e))}$ \\
 \hline
 $X$ \T \B & $6d+\frac{6 a b c}{a b+a c+b c}$  & $\frac{\elg}{9}+\frac{11d}{9}-\frac{2abc}{9 (b c + a (b + c))}$ \\
 \hline
 $XI$ \T \B & $6(c+d)+\frac{8 e f}{e+f}$  & $\frac{\elg}{9}+\frac{11(c+d)}{9}+\frac{2ef}{3(e+f)}$ \\
 \hline
 \small{$XII$} \T \B & $6e+\frac{6 a b c + 6 a b d + 8 a c d + 8 b c d}{a b + (a + b) (c + d)}$   & $\frac{\elg}{9}+\frac{11e}{9}+\frac{6 (a + b) c d - 2 a b (c + d)}{9 (a b + (a + b) (c + d))}$ \\
 \hline
 \hline
 \hline
\end{tabular}
\end{center}  \caption{Values of $\tcg$ and $\vg$ for pm-graphs with $g(\ga)=3$ and $\gc=3$.} \label{tab caseIVb}
\end{table}

\begin{table}
\begin{center}
\begin{tabular}{|c|c|c|}
  \hline
  \T \B & $\lag$ & $\ed$ \\
 \hline
 $I$ \T \B & $\frac{3\elg}{28}$  & $\frac{2\elg}{9}$ \\
 \hline
 $II$ \T \B & $\frac{3\elg}{28}$  & $\frac{2\elg}{9}+\frac{4 a b c d}{9 (b c d + a (c d + b (c + d)))}$ \\
 \hline
 $III$ \T \B & $\frac{3\elg}{28}+\frac{ a b c}{28 (b c + a (b + c))}$  & $\frac{2\elg}{9}+\frac{5a b c}{9 (b c + a (b + c))}$ \\
 \hline
 $IV$ \T \B & $\frac{3\elg}{28}+\frac{c d}{7(c+d)}$  & $\frac{2\elg}{9}+\frac{4c d}{3(c+d)}$ \\
 \hline
$V$ \T \B &  $\frac{3\elg}{28}+\frac{5d}{28}$ & $\frac{2\elg}{9}+\frac{13d}{9}$ \\
 \hline
 $VI$ \T \B & $\frac{3\elg}{28}+\frac{5(d+e)}{28}$  & $\frac{2\elg}{9}+\frac{13(d+e)}{9}$ \\
 \hline
$VII$ \T \B & $\frac{3\elg}{28}+\frac{5c}{28}+\frac{de}{7(d+e)}$  & $\frac{2\elg}{9}+\frac{13c}{9}+\frac{4de}{3(d+e)}$ \\
 \hline
\small{$VIII$} \T \B & $\frac{3\elg}{28}+\frac{a ((c+b) d e + b c (d + e))}{28 (c d e + a (b + c) (d + e) + b (d e + c (d + e)))}$  & $\frac{2\elg}{9}+\frac{4 b c d e + 5 a (c d e + b (d e + c (d + e)))}{9 (c d e + a (b + c) (d + e) + b (d e + c (d + e)))}$ \\
 \hline
 $IX$ \T \B & $\frac{3\elg}{28}+\frac{(b + c) d e + 4 b c (d + e)}{28 (d e + (b + c) (d + e))}$  & $\frac{2\elg}{9}+\frac{5 (b + c) d e + 12 b c (d + e)}{9 (d e + (b + c) (d + e))}$ \\
 \hline
 $X$ \T \B & $\frac{3\elg}{28}+\frac{5d}{28}+\frac{abc}{28 (b c + a (b + c))}$  & $\frac{2\elg}{9}+\frac{13d}{9}+\frac{5abc}{9 (b c + a (b + c))}$ \\
 \hline
 $XI$ \T \B & $\frac{3\elg}{28}+\frac{5(c+d)}{28}+\frac{ef}{7(e+f)}$  & $\frac{2\elg}{9}+\frac{13(c+d)}{9}+\frac{4ef}{3(e+f)}$ \\
 \hline
 \small{$XII$} \T \B & $\frac{3\elg}{28}+\frac{5e}{28}+\frac{4 (a + b) c d + a b (c + d)}{28 (a b + (a + b) (c + d))}$   & $\frac{2\elg}{9}+\frac{13e}{9}+\frac{12 (a + b) c d + 5 a b (c + d)}{9 (a b + (a + b) (c + d))}$ \\
 \hline
\end{tabular}
\end{center}  \caption{Values of $\lag$ and $\ed$ for pm-graphs with $g(\ga)=3$ and $\gc=3$.} \label{tab caseIVc}
\end{table}

\textbf{Pm-graphs of type $XIII$:}

Let $\ga$ be a pm-graph as illustrated in $XIII$ in \figref{fig caseIV}. In this case, we have
$\elg=a+b+c+d+e+f$, $\delta_0(\ga)=\elg$, and $\delta_i(\ga)=0$ for all $i \geq 1$. Moreover,

\begin{align*}
\tg&=\frac{\elg}{12}-\frac{A+2C}{6 D},&
\tcg&=\frac{6A+8B+8C}{D}, \\
\vg&=\frac{\elg}{9}-\frac{2A-6B+7C}{9 D},&
\lag&=\frac{3\elg}{28}+\frac{A+4B}{28 D},
\end{align*}
and
$$ \ed =\frac{2}{9}\elg+\frac{5A+12B+4C}{9 D},$$
where $A=a c d e + b c d e + a c d f + b c d f + a c e f + b c e f + a d e f +
 b d e f$, $B=a b c e + a b d e + a b c f + a b d f$, $C=c d e f$ and
$D=(a + b) c e + (a + b) d e + c d e + (a + b) c f + (a + b) d f + c d f + c e f + d e f)$.

Let $H=c^2 d e+c d^2 e+c d e^2+c^2 d f+c d^2 f+c^2 e f-12 c d e f+d^2 e f+c e^2 f+d e^2 f+c d f^2+c e f^2+d e f^2$,
$N=14 (a + b) c e + 3 c^2 e + 14 (a + b) d e + 3 d^2 e + 3 c e^2 +
 3 d e^2 + 14 (a + b) c f + 3 c^2 f + 14 (a + b) d f + 3 d^2 f +
 3 c f^2 + 3 d f^2$, $M=(c - d)^2 e + (c - d)^2 f + c (e - f)^2 + d (e - f)^2$. We see that $D$, $N$ and $M$ are nonnegative, and
note that $H \geq 0$ by applying Arithmetic-Harmonic Mean inequality for $c$, $d$, $e$ and $f$.

Now we note that
$\vg=\frac{1}{16}\elg+\frac{(a+b)(N+11 M)+14 H+192 a b (c + d) (e + f)}{288 D}$, which implies that $\vg \geq \frac{1}{16}\elg$. This is the sharp lower bound because
$\vg = \frac{1}{16}\elg$ whenever $a=b=0$ and $c=d=e=f$.

\textbf{Pm-graphs of type $XIV$:}

Let $\ga$ be a pm-graph as illustrated in $XIV$ in \figref{fig caseIV}. In this case, we have
$\elg=a+b+c+d+e+f$, $\delta_0(\ga)=\elg$, and $\delta_i(\ga)=0$ for all $i \geq 1$. Moreover,

\begin{align*}
\tg&=\frac{\elg}{12}-\frac{A+2B}{6 C},&
\tcg&=\frac{6A+8B}{C}, \\
\vg&=\frac{\elg}{9}-\frac{2A+7B}{9 C},&
\lag&=\frac{3\elg}{28}+\frac{A}{28 C},
\\  \ed &=\frac{2}{9}\elg+\frac{5A+4B}{9 C}, &  &
\end{align*}

where $A=a b c d + a b c e + a b d e + a c d e + a b c f + a b d f + b c d f +
 a c e f + b c e f + a d e f + b d e f + c d e f$, $B=b c d e + a c d f + a b e f$, and
$C=a b d + a c d + b c d + a b e + a c e + b c e + b d e + c d e +
 a b f + a c f + b c f + a d f + c d f + a e f + b e f + d e f$.

We first show that  $\tg \geq \frac{5}{96} \elg $, where the equality holds whenever $a=b=c=d=e=f$.
We have $\tg=\frac{5}{96} \elg + \frac{3M-7A-20B}{96 C}$
where
$M=a^2 b d+a^2 b e+a^2 b f+a^2 c d+a^2 c e+a^2 c f+a^2 d f+a^2 e f+a b^2 d+a b^2 e+a b^2 f+a b d^2+a b e^2+a b f^2+a c^2 d+a c^2 e+a c^2 f+a c d^2+a c e^2+a c f^2+a d^2 f+a d f^2+a e^2 f+a e f^2+b^2 c d+b^2 c e+b^2 c f+b^2 d e+b^2 e f+b c^2 d+b c^2 e+b c^2 f+b c d^2+b c e^2+b c f^2+b d^2 e+b d e^2+b e^2 f+b e f^2+c^2 d e+c^2 d f+c d^2 e+c d^2 f+c d e^2+c d f^2+d^2 e f+d e^2 f+d e f^2$.

Thus, we see that proving $S:=3M-7A-20B \geq 0$ gives  $\tg \geq \frac{5}{96} \elg $. Now, we have the following tricky equality
\begin{equation*}\label{eqn tau1}
\begin{split}
S&= 2 \big[ b e \left((a-f)^2+(d-c)^2\right)+c d \left((a-f)^2+(b-e)^2\right)+a f \left((e-b)^2+(d-c)^2\right) \big]\\
& \quad + \frac{3}{2}  \big[ b d \left((a-c)^2+(a-e)^2+(c-e)^2\right)+c e \left((b-a)^2+(d-a)^2+(b-d)^2\right)\\
& \quad + a d \left((b-c)^2+(b-f)^2+(c-f)^2\right)+c f \left((a-b)^2+(a-d)^2+(d-b)^2\right)\\
& \quad +b f \left((c-a)^2+(e-a)^2+(c-e)^2\right)+a e \left((b-c)^2+(f-b)^2+(f-c)^2\right)\\
& \quad +e f \left((a-b)^2+(a-d)^2+(d-b)^2\right)+a b \left((d-e)^2+(d-f)^2+(f-e)^2\right)\\
& \quad +a c \left((d-e)^2+(d-f)^2+(e-f)^2\right)+d f \left((a-c)^2+(e-a)^2+(e-c)^2\right)\\
& \quad +d e \left((b-c)^2+(b-f)^2+(c-f)^2\right)+b c \left((d-e)^2+(d-f)^2+(e-f)^2\right)\big]\\
& \quad + \frac{1}{2} \big[ c d (a-b)^2+b e (a-c)^2+a f (b-c)^2+b e (a-d)^2+a f (b-d)^2+c d (a-e)^2\\
& \quad +a f (c-e)^2+a f (d-e)^2+b e \left((c-f)^2+(d-f)^2\right)+c d (b-f)^2+c d (e-f)^2 \big].
\end{split}
\end{equation*}
Thus, $S$ is a sum of positive terms. This gives
\begin{equation}\label{eqn tau0}
\begin{split}
3M-7A-20B \geq 0, \quad  \text{ and so } \quad  \tg \geq \frac{5}{96} \elg.
\end{split}
\end{equation}
Now, we consider $\vg$.
We note that $\vg = \frac{17}{288} \elg$ whenever $\ga$ has equal edge lengths, i.e., if $a=b=c=d=e=f$. Next, we show that this is the sharp lower bound for $\ga$ and so for all pm-graphs of $\gc=3$.

\textbf{Claim:} $\vg \geq \frac{17}{288} \elg$.

\textbf{Proof of Claim:}
It is worth mentioning that we are unable to give a proof of this inequality neither by utilizing arithmetic harmonic mean inequalities partially or fully as in the previous cases nor by using any other well-known inequality in literature. Instead we found the following highly tricky and technical proof after spending extensive time on this problem.

Since $\vg = \frac{17}{288} \elg + \frac{15D-19A-164B}{C}$, where $A$, $B$ and $C$ are as above and $D$ is as follows:

$D=a^2 b d+a^2 b e+a^2 b f+a^2 c d+a^2 c e+a^2 c f+a^2 d f+a^2 e f+a b^2 d+a b^2 e+a b^2 f+a b d^2+a b e^2+a b f^2+a c^2 d+a c^2 e+a c^2 f+a c d^2+a c e^2+a c f^2+a d^2 f+a d f^2+a e^2 f+a e f^2+b^2 c d+b^2 c e+b^2 c f+b^2 d e+b^2 e f+b c^2 d+b c^2 e+b c^2 f+b c d^2+b c e^2+b c f^2+b d^2 e+b d e^2+b e^2 f+b e f^2+c^2 d e+c^2 d f+c d^2 e+c d^2 f+c d e^2+c d f^2+d^2 e f+d e^2 f+d e f^2$

Therefore, it is enough to show that the following inequality holds to prove the claim:
\begin{equation}\label{eqn main inequality0}
\begin{split}
15D-19A-164B \geq 0
\end{split}
\end{equation}
The proof of this inequality consists of eight similar cases that depends on the comparison of the involved variables. The idea is to express $15D-19A-164B$ as sums squares and nonnegative terms.
Lets denote this term by $R$, i.e., we set $R:=15D-19A-164B$.

\textbf{Case I:} Suppose $a \geq f$, $b \geq e$ and $c \geq d$:

We have
\begin{equation*}\label{eqn term1}
\begin{split}
R & =2 \big[ c d (b + e - a - f)^2 + b e (c + d - a - f)^2 + a f (c + d - b - e)^2 \big]\\
& \quad +13 \big[ b e (a-c+d-f)^2+c d (a-b+e-f)^2+a f (-b+c-d+e)^2 \big]\\
& \quad +15 \big[ d e (b - c)^2 + b d (c - e)^2 + d f (c - a)^2 + d a (c - f)^2 + a b (e - f)^2 + a c (d - f)^2\\
& \quad + a e (b - f)^2 + b c (d - e)^2 + b f (a - e)^2 + c e (b - d)^2 + c f (a - d)^2 + e f (a - b)^2 \big]\\
& \quad +11 \big[ c d (a-f) (b-e)+b e (a-f) (c-d)+a f (b-e) (c-d) \big]\\
&\quad +15 T_1,
\end{split}
\end{equation*}
where
$T_1= a^2 b d+a^2 c e+a b^2 d+a b d^2-2 a b d e-2 a b d f+a c^2 e-2 a c d e+a c e^2-2 a c e f+b^2 c f+b c^2 f-2 b c d f-2 b c e f+b c f^2+d^2 e f+d e^2 f+d e f^2$.

By the assumptions in this case, we have $a-f \geq 0$, $b-e \geq 0$ and $c-d \geq 0$. Therefore, to prove $R \geq 0$, it will be enough to show $T_1 \geq 0$. Again by the assumptions, we can write $a=f+k$, $b=e+m$ and $c=d+n$ for some nonnegative real numbers $k$, $m$ and $n$. Now substituting these into $T_1$ gives

$T_1=d^2 k m+2 d e k^2+2 d f m^2+d k^2 m+d k m^2+e^2 k n+2 e f n^2+e k^2 n+e k n^2+f^2 m n+f m^2 n+f m n^2$,
which clearly shows that $T_1 \geq 0$. Hence, $R \geq 0$ in this case.

\textbf{Case II:} Suppose $ a \geq f$, $ b \geq e $ and $d \geq c$:

We have
\begin{equation*}\label{eqn term2}
\begin{split}
R & =2 \big[ a f (-b+c+d-e)^2+b e (-a+c+d-f)^2+c d (-a+b+e-f)^2 \big]\\
& \quad +13 \big[ b e (a+c-d-f)^2+c d (a-b+e-f)^2+a f (-b-c+d+e)^2 \big]\\
& \quad +15 \big[ a e (b-f)^2+a b (e-f)^2+b f (a-e)^2+e f (a-b)^2+d f (c-a)^2+a d (c-f)^2\\
& \quad +a c (d-f)^2+c f (a-d)^2+b d (c-e)^2+b c (d-e)^2+c e (b-d)^2+d e (b-c)^2  \big]\\
& \quad +11 \big[ c d (a-f) (b-e)+b e (a-f) (d-c)+a f (b-e) (d-c) \big]\\
&\quad +15 T_2,
\end{split}
\end{equation*}
where
$T_2= a^2 b d+a^2 c e+a b^2 d-2 a b c e-2 a b c f+a b d^2+a c^2 e-2 a c d e+a c e^2-2 a d e f+b^2 c f+b c^2 f-2 b c d f+b c f^2-2 b d e f+d^2 e f+d e^2 f+d e f^2$.

By the assumptions in this case, we have $a-f \geq 0$, $b-e \geq 0$ and $d-c \geq 0$. Therefore, to prove $R \geq 0$, it will be enough to show $T_2 \geq 0$. Again by the assumptions, we can write $a=f+k$, $b=e+m$ and $d=c+n$ for some nonnegative real numbers $k$, $m$ and $n$. Now substituting these into $T_2$ gives
$T_2=c^2 k m+2 c e k^2+2 c f m^2+c k^2 m+c k m^2+2 c k m n+e^2 k n+2 e f n^2+e k^2 n+2 e k m n+e k n^2+f^2 m n+2 f k m n+f m^2 n+f m n^2+k^2 m n+k m^2 n+k m n^2$. This shows that $T_2 \geq 0$. Hence, $R \geq 0$ in this case.

\textbf{Case III:} Suppose $ a\geq f$, $e \geq b$ and $ c \geq d $:

We have
\begin{equation*}\label{eqn term3}
\begin{split}
R & =2 \big[ a f (-b+c+d-e)^2+b e (-a+c+d-f)^2+c d (-a+b+e-f)^2 \big]\\
& \quad +13 \big[ c d (a - f + b - e)^2 + b e (a - f - c + d)^2 + a f (c - d + b - e)^2 \big]\\
& \quad +15 \big[ a e (b-f)^2+a b (e-f)^2+b f (a-e)^2+e f (a-b)^2+d f (c-a)^2+a d (c-f)^2\\
& \quad +a c (d-f)^2+c f (a-d)^2+b d (c-e)^2+b c (d-e)^2+c e (b-d)^2+d e (b-c)^2  \big]\\
& \quad +11 \big[ c d (a-f) (e-b)+b e (a-f) (c-d)+a f (e-b) (c-d) \big]\\
&\quad +15 T_3,
\end{split}
\end{equation*}
where
$T_3= a^2 b d+a^2 c e+a b^2 d-2 a b c d-2 a b c f+a b d^2-2 a b d e+a c^2 e+a c e^2-2 a d e f+b^2 c f+b c^2 f-2 b c e f+b c f^2-2 c d e f+d^2 e f+d e^2 f+d e f^2$.

By the assumptions in this case, we have $ a-f \geq 0$, $ e-b \geq 0$ and $c-d \geq 0$. Therefore, to prove $R \geq 0$, it will be enough to show $T_3 \geq 0$. Again by the assumptions, we can write $a=f+k$, $e=b+m$ and $c=d+n$ for some nonnegative real numbers $k$, $m$ and $n$. Now substituting these into $T_3$ gives
$T_3=b^2 k n+2 b d k^2+2 b f n^2+b k^2 n+2 b k m n+b k n^2+d^2 k m+2 d f m^2+d k^2 m+d k m^2+2 d k m n+f^2 m n+2 f k m n+f m^2 n+f m n^2+k^2 m n+k m^2 n+k m n^2$,
which shows that $T_3 \geq 0$. Hence, $R \geq 0$ in this case.

\textbf{Case IV:} Suppose $a \geq f$, $e \geq b$ and $d \geq c$:

We have
\begin{equation*}\label{eqn term4}
\begin{split}
R & =2 \big[ a f (-b+c+d-e)^2+b e (-a+c+d-f)^2+c d (-a+b+e-f)^2 \big]\\
& \quad +13 \big[ b e (a+c-d-f)^2+c d (a+b-e-f)^2+a f (b-c+d-e)^2 \big]\\
& \quad +15 \big[ a e (b-f)^2+a b (e-f)^2+b f (a-e)^2+e f (a-b)^2+d f (c-a)^2+a d (c-f)^2\\
& \quad +a c (d-f)^2+c f (a-d)^2+b d (c-e)^2+b c (d-e)^2+c e (b-d)^2+d e (b-c)^2  \big]\\
& \quad +11 \big[ c d (a-f) (e-b)+b e (a-f) (d-c)+a f (e-b) (d-c) \big]\\
&\quad +15 T_4,
\end{split}
\end{equation*}
where
$T_4= a^2 b d+a^2 c e+a b^2 d-2 a b c d-2 a b c e+a b d^2-2 a b d f+a c^2 e+a c e^2-2 a c e f+b^2 c f+b c^2 f+b c f^2-2 b d e f-2 c d e f+d^2 e f+d e^2 f+d e f^2$.

By the assumptions in this case, we have $a-f \geq 0$, $e-b \geq 0$ and $d-c \geq 0$. Therefore, to prove $R \geq 0$, it will be enough to show $T_4 \geq 0$. Again by the assumptions, we can write $a=f+k$, $e=b+m$ and $d=c+n$ for some nonnegative real numbers $k$, $m$ and $n$. Now substituting these into $T_4$ gives
$T_4=b^2 k n+2 b c k^2+2 b f n^2+b k^2 n+b k n^2+c^2 k m+2 c f m^2+c k^2 m+c k m^2+f^2 m n+f m^2 n+f m n^2$,
which clearly shows that $T_4 \geq 0$. Hence, $R \geq 0$ in this case.

\textbf{Case V:} Suppose $ f \geq a$, $b \geq e$ and $c \geq d$:

We have
\begin{equation*}\label{eqn term5}
\begin{split}
R & =2 \big[ a f (-b+c+d-e)^2+b e (-a+c+d-f)^2+c d (-a+b+e-f)^2  \big]\\
& \quad +13 \big[ b e (-a-c+d+f)^2+c d (-a-b+e+f)^2+a f (-b+c-d+e)^2 \big]\\
& \quad +15 \big[ a e (b-f)^2+a b (e-f)^2+b f (a-e)^2+e f (a-b)^2+d f (c-a)^2+a d (c-f)^2 \\
& \quad +a c (d-f)^2+c f (a-d)^2+b d (c-e)^2+b c (d-e)^2+c e (b-d)^2+d e (b-c)^2  \big]\\
& \quad +11 \big[ a f (b-e) (c-d)+c d (f-a) (b-e)+b e (f-a) (c-d) \big]\\
&\quad +15 T_5,
\end{split}
\end{equation*}
where
$T_5= a^2 b d+a^2 c e+a b^2 d-2 a b c d-2 a b c e+a b d^2-2 a b d f+a c^2 e+a c e^2-2 a c e f+b^2 c f+b c^2 f+b c f^2-2 b d e f-2 c d e f+d^2 e f+d e^2 f+d e f^2$.

By the assumptions in this case, we have $ f-a \geq 0$, $b-e \geq 0$ and $c-d \geq 0$. Therefore, to prove $R \geq 0$, it will be enough to show $T_5 \geq 0$. Again by the assumptions, we can write $f=a+k$, $b=e+m$ and $c=d+n$ for some nonnegative real numbers $k$, $m$ and $n$. Now substituting these into $T_5$ gives
$T_5=a^2 m n+2 a d m^2+2 a e n^2+2 a k m n+a m^2 n+a m n^2+d^2 k m+2 d e k^2+d k^2 m+d k m^2+2 d k m n+e^2 k n+e k^2 n+2 e k m n+e k n^2+k^2 m n+k m^2 n+k m n^2$. Thus, $T_5 \geq 0$. Hence, $R \geq 0$ in this case.

\textbf{Case VI:} Suppose $ f \geq a $, $b \geq e$ and $d \geq c$:

We have

\begin{equation*}\label{eqn term6}
\begin{split}
R & =2 \big[ a f (-b+c+d-e)^2+b e (-a+c+d-f)^2+c d (-a+b+e-f)^2 \big]\\
& \quad +13 \big[ b e (-a+c-d+f)^2+c d (-a-b+e+f)^2+a f (-b-c+d+e)^2 \big]\\
& \quad +15 \big[ a e (b-f)^2+a b (e-f)^2+b f (a-e)^2+e f (a-b)^2+d f (c-a)^2+a d (c-f)^2\\
& \quad +a c (d-f)^2+c f (a-d)^2+b d (c-e)^2+b c (d-e)^2+c e (b-d)^2+d e (b-c)^2  \big]\\
& \quad +11 \big[ a f (b-e) (d-c)+c d (f-a) (b-e)+b e (f-a) (d-c)  \big]\\
&\quad +15 T_6,
\end{split}
\end{equation*}
where
$T_6= a^2 b d+a^2 c e+a b^2 d-2 a b c d-2 a b c f+a b d^2-2 a b d e+a c^2 e+a c e^2-2 a d e f+b^2 c f+b c^2 f-2 b c e f+b c f^2-2 c d e f+d^2 e f+d e^2 f+d e f^2$.

By the assumptions in this case, we have $f-a \geq 0$, $b-e \geq 0$ and $d-c \geq 0$. Therefore, to prove $R \geq 0$, it will be enough to show $T_6 \geq 0$. Again by the assumptions, we can write $f=a+k$, $b=e+m$ and $d=c+n$ for some nonnegative real numbers $k$, $m$ and $n$. Now substituting these into $T_6$ gives
$T_6=a^2 m n+2 a c m^2+2 a e n^2+a m^2 n+a m n^2+c^2 k m+2 c e k^2+c k^2 m+c k m^2+e^2 k n+e k^2 n+e k n^2$,
so $T_6 \geq 0$. Hence, $R \geq 0$ in this case.

\textbf{Case VII:} Suppose $ f \geq a $, $e \geq b $ and $c \geq d$:

We have
\begin{equation*}\label{eqn term7}
\begin{split}
R & =2 \big[ a f (-b+c+d-e)^2+b e (-a+c+d-f)^2+c d (-a+b+e-f)^2 \big]\\
& \quad +13 \big[ b e (-a-c+d+f)^2+c d (-a+b-e+f)^2+a f (b+c-d-e)^2 \big]\\
& \quad +15 \big[ a e (b-f)^2+a b (e-f)^2+b f (a-e)^2+e f (a-b)^2+d f (c-a)^2+a d (c-f)^2\\
& \quad +a c (d-f)^2+c f (a-d)^2+b d (c-e)^2+b c (d-e)^2+c e (b-d)^2+d e (b-c)^2  \big]\\
& \quad +11 \big[ a f (e-b) (c-d)+b e (f-a) (c-d)+c d (f-a) (e-b) \big]\\
&\quad +15 T_7,
\end{split}
\end{equation*}
where
$T_7=a^2 b d+a^2 c e+a b^2 d-2 a b c e-2 a b c f+a b d^2+a c^2 e-2 a c d e+a c e^2-2 a d e f+b^2 c f+b c^2 f-2 b c d f+b c f^2-2 b d e f+d^2 e f+d e^2 f+d e f^2 $

By the assumptions in this case, we have $ f-a \geq 0$, $e-b \geq 0$ and $c-d \geq 0$. Therefore, to prove $R \geq 0$, it will be enough to show $T_7 \geq 0$. Again by the assumptions, we can write $f=a+k$, $e=b+m$ and $c=d+n$ for some nonnegative real numbers $k$, $m$ and $n$. Now substituting these into $T_7$ gives
$T_7=a^2 m n+2 a b n^2+2 a d m^2+a m^2 n+a m n^2+b^2 k n+2 b d k^2+b k^2 n+b k n^2+d^2 k m+d k^2 m+d k m^2$,
which clearly shows that $T_7 \geq 0$. Hence, $R \geq 0$ in this case.

\textbf{Case VIII:} Suppose $f \geq a$, $e \geq b$ and $d \geq c$:

We have
\begin{equation*}\label{eqn term8}
\begin{split}
R & =2 \big[ a f (-b+c+d-e)^2+b e (-a+c+d-f)^2+c d (-a+b+e-f)^2 \big]\\
& \quad +13 \big[ b e (-a+c-d+f)^2+c d (-a+b-e+f)^2+a f (b-c+d-e)^2 \big]\\
& \quad +15 \big[ a e (b-f)^2+a b (e-f)^2+b f (a-e)^2+e f (a-b)^2+d f (c-a)^2+a d (c-f)^2\\
& \quad +a c (d-f)^2+c f (a-d)^2+b d (c-e)^2+b c (d-e)^2+c e (b-d)^2+d e (b-c)^2  \big]\\
& \quad +11 \big[ a f (e-b) (d-c)+b e (f-a) (d-c)+c d (f-a) (e-b) \big]\\
&\quad +15 T_8,
\end{split}
\end{equation*}
where
$T_8=a^2 b d+a^2 c e+a b^2 d+a b d^2-2 a b d e-2 a b d f+a c^2 e-2 a c d e+a c e^2-2 a c e f+b^2 c f+b c^2 f-2 b c d f-2 b c e f+b c f^2+d^2 e f+d e^2 f+d e f^2$.

By the assumptions in this case, we have $f-a \geq 0$, $e-b \geq 0$ and $d-c \geq 0$. Therefore, to prove $R \geq 0$, it will be enough to show $T_8 \geq 0$. Again by the assumptions, we can write $f=a+k$, $e=b+m$ and $d=c+n$ for some nonnegative real numbers $k$, $m$ and $n$. Now substituting these into $T_8$ gives
$T_8= a^2 m n+2 a b n^2+2 a c m^2+2 a k m n+a m^2 n+a m n^2+b^2 k n+2 b c k^2+b k^2 n+2 b k m n+b k n^2+c^2 k m+c k^2 m+c k m^2+2 c k m n+k^2 m n+k m^2 n+k m n^2$. This clearly shows that $T_8 \geq 0$. Hence, $R \geq 0$ in this case.





Next, we give a summary of the inequalities that we established so far.
If $(\ga,\bq)$ is a pm-graph of $\gc=3$ that is not a single point, then we showed that we have the following equalities and sharp lower bounds:

We have $\lag =\frac{2}{7} \elg$ and $\ed =\frac{5}{3} \elg$ if $g=0$, and $\lag \geq \frac{3}{28} \elg$ and $\ed \geq \frac{2}{9} \elg$ if $1 \leq g \leq 3$.

If $g=0$, $\vg =\frac{4}{3} \elg$, and if $g=1$, $\vg  \geq \frac{1}{9} \elg$. When $g=2$, $\vg  \geq \frac{7}{81} \elg$. Finally,
$\vg  \geq \frac{17}{288} \elg$ if $g=3$.

For any nonnegative six real numbers $a$, $b$, $c$, $d$, $e$ and $f$, we showed (for pm-graph of type XIV in  \figref{fig caseIV})
that
\begin{equation}\label{eqn main inequality}
\begin{split}
\frac{15}{32}(a+b+c+d+e+f)-\frac{7(a b e f+a c d f+b c d e)+2 A}{C} \geq 0,
\end{split}
\end{equation}
where
$A=c d(b+e)(a+f)+b e(c+d)(a+f)+a f(c+d)(b+e)$
and
$C=c d (a+b+e+f)+a f (b+c+d+e)+b e (a+c+d+f)+a b d+a c e+b c f+d e f$.

The equality in (\ref{eqn main inequality}) holds if $a=b=c=d=e=f>0$. Similarly, we can rewrite (\ref{eqn tau0}) as follows
\begin{equation}\label{eqn main inequality3}
\begin{split}
\frac{1}{32}(a+b+c+d+e+f)-\frac{2(a b e f+a c d f+b c d e)+ A}{C} \geq 0,
\end{split}
\end{equation}
Moreover, if $a+b+c+d+e+f=1$, we can rewrite inequalities (\ref{eqn main inequality}) and (\ref{eqn main inequality3}) as follows:
\begin{equation}\label{eqn main inequality2}
\begin{split}
&15 C -224A-64(a b e f+a c d f+b c d e) \geq 0,\\
&3C-16A-32(a b e f+a c d f+b c d e) \geq 0,
\end{split}
\end{equation}
where the equalities holds iff $a=b=c=d=e=f=\frac{1}{6}.$

\textbf{Acknowledgements:} This work is supported by The Scientific and Technological Research Council of Turkey-TUBITAK (Project No: 110T686).


\end{document}